\documentclass[12pt]{article}
\usepackage{amsxtra,amssymb,amsthm,amsmath,latexsym}

\theoremstyle{plain}

\def\R{{\mathbb R}}

\def\oH{{\overset{\circ}{H}}}
\def\oH1{{\overset{\circ}{H}\kern-.02in{}^1}}

\def\bee{\begin{equation*}}
\def\eee{\end{equation*}}
\def\be{\begin{equation}}
\def\ee{\end{equation}}

\begin{document}

Published in: Far East Journ. of Appl. Math.,  116, N3, (2023), pp. 215-227.

\title{ Boundary values of analytic functions
}

\author{Alexander G. Ramm\\
 Department  of Mathematics, Kansas State University, \\
 Manhattan, KS 66506, USA\\
ramm@ksu.edu\\
http://www.math.ksu.edu/\,$\sim$\,ramm}

\date{}
\maketitle\thispagestyle{empty}

\begin{abstract}

Let $D$ be a connected bounded domain in $\R^2$, $S$ be its boundary
which is closed, connected and smooth. Let $\Phi(z)=\frac 1 {2\pi i}\int_S\frac{f(s)ds}{s-z}$, $f\in L^1(S)$, $z=x+iy$.
Boundary values of $\Phi(z)$ on $S$ are studied. The function $\Phi(t)$,
$t\in S$, is defined in a new way.  Necessary and sufficient conditions are
given for $f\in L^1(S)$ to be boundary value of an analytic in $D$ function.
The Sokhotsky-Plemelj formulas are derived for  $f\in L^1(S)$. 
\end{abstract}

\section{Introduction}\label{S:1}
 Let $D$ be a connected bounded domain on the complex plane, $S$ be its  boundary, which is closed and $C^{1,a}$-smooth, $0<a\le 1$.
Consider an analytic function in $D$, defined as
\be\label{e1}  
\Phi(z)=c\int_S\frac{f(s)}{s-z}ds, \quad c:=\frac 1 {2\pi i}.
\ee
We assume that $f\in L^1(S)$. This is {\em the basic new assumption}: in the literature it was assumed that $f$ is H$\ddot{o}$lder-continuous or
$f\in L^p(S),\,\, p>1$, see \cite{G}, \cite{MP}. In \cite{C} there is a result for $f\in L^1(S)$, which we mention in Remark 3.

It is of great interest to have a proof
of the Sokhotsky-Plemelj formulas for $f\in L^1(S)$ and of the the relation
\eqref{e8}  for $f\in L^1(S)$, see below. The contents of Chapter 2 in
 \cite{G}  is based on the relation \eqref{e8},
which is proved in \cite{G} under the assumption that $f$ is H$\ddot{o}$lder-continuous, $f\in H^\mu(S)$, $\mu\in (0, 1]$.  We prove that this relation holds almost everywhere (a.e.) in the sense of the Lebesgue measure on $S$.
The Sokhotsky-Plemelj formulas for $f\in L^1(S)$ are important in applications to singular integral operators and boundary value problems,
\cite{G}, \cite{M}.

By $D^+$ we denote $D$, by $D^-=D'$
we denote $\R^2\setminus \bar{D}$, the $\bar{D}$ is the closure of $D$.
 The function $\Phi$ is analytic in $D$
and in $D'$, $\Phi(\infty)=0$.  Boundary value of $\Phi$ on $S$ has been
studied  by many authors, \cite{C}--\cite{P}.  In \cite{G} and 
\cite{M} it is assumed that $f\in H^\mu(S)$. In \cite{MP}
and in \cite{Kh} it is assumed that $f\in L^p(S)$, $p>1$.  The case $p=1$ is not discussed in these books. In \cite{C} there are results related to $p=1$
related to the Cauchy principal value definition of $\Phi(t)$. We give a new
definition of the singular integral $\Phi(t)$, see \eqref{e9b}.

Define $\psi(z)$ by formula \eqref{e2}.
Denote the limiting values of $\Phi(z)$ and $\psi(z)$, when $z\to t\in S$ along the non-tangential to $S$ directions, by 
$\Phi^+(t)$ and $\psi^+(t)$
when $z\in D, \, z\to t$, and by  $\Phi^-(t)$ and $\psi^-(t)$
when $z\in D', \, z\to t$. By $\Phi(t)$ and $\psi(t)$ we denote the values of
$\Phi(z)$ and $\psi(z)$ when $z=t\in S$. 
Let us assume that $z\to t$ along the unit normal $N_t$ to $S$, where the normal is directed out of $D$. 

  We are interested in the following question: 

 {\em How to 
characterize  functions $f\in L^1(S)$ which are boundary values of
analytic functions in $D$? }

We answer this question in Theorem 1 (see below).  

Rewrite \eqref{e1} as
\be\label{e2}  
\Phi(z)=f(t)\nu(z) + \psi(z),\quad \psi(z):=
c\int_S\frac{f(s)-f(t)}{s-z}ds, \quad t\in S,\quad \nu(z):=c\int_S\frac{ds}{s-z},
\ee
and 
\begin{equation}\label{e3}
\nu(z)=\begin{cases}
1,\quad z\in D,\\
\frac 1 2, \quad z\in S,\\
0, \quad z\in D'.
\end{cases}
\end{equation}
One has
\be\label{e4}  
\Phi^+(t)=\lim_{z\to t, z\in D} \Phi(z)=f(t)+\psi^+(t)
\ee
where $\psi^{+}(t)=\lim_{z\to t, z\in D} \psi(z)$ and
\be\label{e5}  
\Phi^-(t)=\lim_{z\to t, z\in D'} \Phi(z)=\psi^-(t).
\ee
If $t\in S$, then one gets (see equation \eqref{e3}, the line $z\in S$) :
\be\label{e6}  
\Phi(t)=\frac {f(t)} 2+\psi(t):=\frac {f(t)} 2+
c\int_S\frac{f(s)-f(t)}{s-t}ds.
\ee

The $\psi(t)$ is the value of $\psi(z)$ at $z=t$.  The $\psi(t)$ and $\Phi(t)$ are understood as in Remark 3, see below.

If some equation holds almost everywhere with respect to the Lebesgue measure on $S$, then we write that this equation holds $a.e.$

From formulas \eqref{e3}--\eqref{e6} one derives:
\be\label{e7}  
\Phi^+(t)-\Phi^-(t)=f(t)+\psi^+(t)-\psi^-(t)\,\, a.e., \quad \Phi^+(t)+\Phi^-(t)=f(t) +\psi^+(t)+\psi^-(t)\,\, a.e.
\ee

In Lemma 1, see below, we prove that $\psi^+(t)=\psi^-(t)=\psi(t)\,\, a.e.$
Therefore, formula \eqref{e7} can be rewritten as:
 \be\label{e7a}  
\Phi^+(t)-\Phi^-(t)=f(t)\,\, a.e., \quad \Phi^+(t)+\Phi^-(t)=f(t) +2\psi(t)\,\, a.e.
\ee

We will need three lemmas. Lemma 1 is proved in Section 2.

{\bf Lemma 1.} {\em If $f\in L^1(S)$ and $S$ is $C^{1,a}-$smooth, 
$0<a\le 1$, then
\be\label{e8}  
\psi^+(t)=\psi^-(t)=\psi(t) \,\,\, a.e.
\ee
}
Let $s,t\in S$

{\bf Lemma 2.} {\em One has
\be\label{e9}  
\frac 1 {s-t\pm i0}=\frac 1 {s-t} \mp i\pi \delta(s-t),
\ee
where $\delta(s-t)$ is the delta-function.
}

{\em Proof.} Formula \eqref{e9} is understood as follows.  Let $f\in L^1(S)$,
$\phi\in H^\mu(S)$, $N_t$ be a unit normal to $S$ directed out of $D$.  Let 
$$\lim_{\epsilon \to +0}\int_S dt \phi(t)\int_S\frac{f(s)ds}{s-t-i\epsilon N_t}:=\int_S dt \phi(t)\int_S\frac{f(s)ds}{s-t-i0}.$$
Then 
\be\label{e9a} 
\int_S dt \phi(t)\int_S\frac{f(s)ds}{s-t-i0}=\int_S dt \phi(t)\int_S\frac{f(s)ds}{s-t}+i\pi \int_S\phi(t) f(t)dt, \quad \forall \phi\in H^\mu(S),
\ee
and
\be\label{e9b} 
\int_S dt \phi(t)\int_S\frac{f(s)ds}{s-t}=\int_S ds f(s) \int_S  \frac{\phi(t)dt}{s-t},  \quad \forall \phi\in H^\mu(S).
\ee
Any $f\in L^1(S)$ is uniquely identified with the linear functional in $L^1(S)$ defined as follows: $(f,\phi)=
\int_Sf(s)\phi(s)ds$, $\forall \phi\in H^\mu(S)$, because the set 
$\{\phi\}|_{\forall \phi\in H^{\mu(S)}}$ is dense in $L^1(S)$.

It is known, see \cite{G}, p. 52, that $B\phi:=\frac 1 {i\pi}\int_S\frac{\phi(s)ds}{s-t}\in H^\mu(S)$, provided that $\phi\in H^\mu(S)$
and $\mu\in (0,1)$. Moreover, the range of $B$ is equal to $H^\mu(S)$
when $\phi$ runs through $H^\mu(S)$, see  \cite{G}, p. 178. This follows from the inversion formulas: if $\phi, \psi\in H^\mu(S)$ and $B\phi=\psi$,
then $\phi=B\psi$, see \cite{M}, p. 115, so $B^2=I$, the identity operator. Therefore, the right side of \eqref{e9b} defines $f\in L^1(S)$ uniquely. The integrals on the left and right sides of \eqref{e9b} are equal to the integral $\int_{S\times S}\frac{f(s) \phi(t)ds dt}{s-t}$, which is absolutely convergent on $S\times S$.
We have proved formula \eqref{e9} with the minus sign. Similarly one proves this formula with the plus sign.
Lemma 2 is proved.  \hfill$\Box$ 

One can interpret formula \eqref{e9a} as follows: if $f\in L^1(S)$, then
\be\label{e9c} 
\int_S\frac{f(s)ds}{s-t-i0}=\int_S \frac {f(s)ds}{s-t} +i\pi f(t) \quad a.e.,
\ee
where the first integral on the right is understood as in formula \eqref{e9b},
see also Remark 3 below, and $f(t)$ is defined $a.e$: almost everywhere  on $S$ with respect to the Lebesgue measure on $S$. 

Usually the delta-function is defined as a linear continuous functional on the space of continuous functions, $\int_S f(s)\delta(s-t)ds=f(t)$, $f\in C(S)$. One may define
delta-function $\delta(s-t)$ on $L^1(S)$ as a linear bounded functional on $L^1(S)$,  $\int_S f(s)\delta(s-t)ds=f(t)\quad a.e.$, $f\in L^(S)$. In place of  $f(t)$ one has to put $f(t) \, a.e.$ because $f\in L^1(S)$ is defined a.e.
with respect to the Lebesgue measure on $S$. By Lusin's theorem, see \cite{S}, p. 157,  if $f\in L^1(S)$, then $f$ is continuous on a subset $S_\epsilon$ of $S$, $|meas S-meas S_\epsilon|<\epsilon$, where $\epsilon>0$ can be chosen arbitrarily small and $meas S$ is the Lebesgue measure of the set $S$. One can consider in this setting the
delta-function $\delta(s-t)$ as a kernel of the identity operator in $L^1(S)$.

In \cite{GS}, p. 83, there is a formula $\frac 1{x-i0}=\frac 1 x +i\pi \delta(x)$
understood in the sense of distributions. The formula in Lemma 2 is  
of the similar type. 

From Lemma 1 and formula \eqref{e7}  one derives the Sokhotsky-Plemelj formulas for $f\in L^1(S)$, see Lemma 3.
 These formulas are derived in \cite{G} and \cite{M} under the assumption that $f$ is H$\ddot{o}$lder-continuous. Under such an assumption, these formulas hold everywhere, not almost everywhere.

{\bf Lemma 3.} {\em If $f\in L^1(S)$, then the Sokhotsky-Plemelj formulas
hold:
\be\label{e10}  
\Phi^+(t)=\Phi(t) + \frac{f(t)} 2\,\, a.e., \quad \Phi^-(t)=\Phi(t) - \frac{f(t)} 2\,\,a.e.
\ee
}

{\bf Proof of Lemma 3.} Formulas \eqref{e10} follow from formulas \eqref{e4}--\eqref{e8}. Indeed, by formulas \eqref{e7} and \eqref{e8}
one concludes that $\Phi^+(t)-\Phi^-(t)=f(t)\, a.e.$ and 
 $\Phi^+(t) +\Phi^-(t)=f(t)+2\psi(t)\, a.e.$ Since $\Phi(t)=\Phi^+(t)-\frac{f(t)}2=\Phi^-(t)+ \frac{f(t)}2$,  formulas \eqref{e10} follow.
Lemma 3 is proved. \hfill$\Box$

Let  $z\in D$. The following question is of interest:

{\em When is the boundary value  $\Phi^+(t)$ of $\Phi (z)$ on $S$ equal to $f$ a.e.?}

From formula \eqref{e4}  the answer follows immediately:

{\em  $\Phi^+(t)=f(t)$ if and only if (iff) $\psi^+(t)=0 \,\,a.e.$}

Equation \eqref{e10} yields another necessary and sufficient condition:

{\em  $\Phi^+(t)=f(t)$ iff  $\Phi(z)=0, \,\,z\in D'$ and $f(t)=\Phi(t)+\frac {f(t)}2$, or
\be\label{e11}  
  f(t)=\frac 1 {i\pi}\int_S\frac{f(s)}{s-t}ds \quad a.e. 
\ee
}
If equation \eqref{e11} holds, then $\Phi^-(t)=0$ and, consequently, 
$ \Phi(z)=0\,\,\text{if}\,\, z\in D'$. 
Equation  \eqref{e11} is a singular integral equation for $f$.
It is clear from this equation that  the singular integral operator  in this equation maps every solution $f\in L^1(S)$ of this equation into itself.

If one wants to formulate a necessary and sufficient condition for $f(s)\in L^1(S)$ to be a boundary value of an analytic in $D'$ function $\Phi(z),\,\,
\Phi(\infty)=0$, then an argument, similar to the above yields the following
conditions:
\be\label{e11a}  
  f(t)=-\frac 1 {i\pi}\int_S\frac{f(s)}{s-t}ds \quad a.e.
\ee  
If equation \eqref{e11a} holds, then $\Phi^+(t)=0$ and, consequently, 
$\Phi(z)=0  \,\, \text{if}\,\, z\in D$.  

{\bf Remark 1.} Integral equation in formula \eqref{e11} has infinitely many linearly independent solutions:
every boundary value of an analytic function $\Phi(z)$ in $D$ such that
$\Phi(z)=0$ in $D'$ solves this integral equation.

{\bf Remark 2.} If $\Phi^+(t)=f(t)\,\, a.e., f\in L^1(S)$,
 then $Bf\in L^1(S)$, where
$Bf:=\frac 1 {i\pi}\int_S\frac{f(s)}{s-t}ds\,\, a.e.$
If $-\Phi^-(t)=f(t)\,\, a.e., f\in L^1(S)$ and $\Phi(\infty)=0$,  then $Bf\in L^1(S)$. Since for some $f\in L^1(S)$ one does not have $\Phi(z)=0, \,\, z\in D'$, it follows that not every $f\in L^1(S)$ is a boundary value of an
analytic function in $D$.

{\bf Remark 3.} Let us define $\Phi(t)$ as follows.

 For any $\phi\in H^\mu(S)$, $\mu\in (0,1)$, we define $\Phi(t)$ by the relation:
\be\label{e11b}  
 \int_S\Phi(t)\phi(t) dt=\int_S f(s) \eta(s) ds, \quad  \forall \phi\in H^\mu(S),
 \ee
  where
$\eta(s):=c\int_S\frac{\phi(t)}{s-t}dt:=B\phi$. 

 If $\phi\in H^\mu(S)$, $\mu\in (0,1)$, then the function $\eta(s)\in H^\mu(S)$, provided that $S$ is  smooth. This is proved in \cite{G}, p. 53.
 
 The singular integral equation  
\be\label{e11c}   
 \frac 1{i\pi}\int_S\frac {\phi(s)ds}{s-t}=h(t), \quad \phi\in H^\mu(S),\,\,
 \mu\in (0,1),
 \ee
 can be solved uniquely and explicitly for any $h\in H^\mu(S)$: 
 \be\label{e11d}   \frac 1{i\pi}\int_S\frac {h(s)ds}{s-t}=\phi(t),
 \ee
 and its solution belongs to $H^\mu(S)$.
  This is proved in \cite{G}, p. 66. 
  
  Consequently, any $h\in H^\mu(S)$, $\mu\in (0,1)$,  can be represented by the integral \eqref{e11c} with $\phi\in H^\mu(S)$. This is in contrast to the 
 integral equation \eqref{e11c} with $\phi\in L^1(S)$: not for every such $\phi$ the $h$ belongs to $L^1(S)$, see Remark 5 below.
 
In \cite{C} it is proved that $\lim_{\epsilon\to 0} \Phi_\epsilon(t)$ exists
for a.e. $t\in S$ if $f\in L^1(S)$, where $\Phi_{\epsilon}(t):=c\int_{|t-s|>\epsilon} \frac{f(s)}{s-t}ds$.

There exists a sequence $f_n\in H^\mu(S), \,\, lim_{n\to \infty} \|f_n-f\|=0$, where the norm is $L^1(S)$-norm.
Clearly $  \int_S\phi(t)\Big(\int_S\frac {f_n(s)ds}{s-t}\Big)dt=\int_Sf_n(s)
\Big(\int_S\frac{\phi(t)dt}{s-t}\Big)ds, \quad \forall \phi\in H^\mu(S).$
The limit of the right side, as $n\to 0$, exists, so the limit of the left side exists and is equal to that of the right side. This leads to the formula \eqref{e11e} below. An alternative derivation of this formula goes as follows. 
Since $\phi(t)f(s)\in L^1(S\times S)$, the operator 
$\int_{S\times S} \phi(t) f(s)|s-t|^{-1}dsdt$ is a bounded operator
in $L^1(S\times S)$, see \cite{GT}, p.159. Consequently, the change of the order of integration in the double integral:
\be\label{e11e}  \int_S\phi(t)\Big(\int_S\frac {f(s)ds}{s-t}\Big)dt=\int_Sf(s)
\Big(\int_S\frac{\phi(t)dt}{s-t}\Big)ds, \quad \forall \phi\in H^\mu(S),
\ee
is legitimate.

Let us summarize the results:

{\bf Theorem 1.} {\em A necessary and sufficient condition for $f\in L^1(S)$
to be the boundary value of an analytic in $D$ function $\Phi(z)$ is
equation \eqref{e11}.

A necessary and sufficient condition for $f\in L^1(S)$
to be the boundary value of an analytic in $D'$ function $\Phi(z), \,\Phi(\infty)=0$, is
equation \eqref{e11a}. 
}

\section{Proof of Lemma 1.}\label{S:2}
 
 Write formula \eqref{e2} for $\psi$ as
 \be\label{e12}  
  \psi^{\pm}(t)= c\int_S\frac{f(s)-f(t)}{s-(t\pm i0)}ds.
\ee
 Using Lemma 2, one obtains:
 \be\label{e13}  
  \psi^{\pm}(t)=c\int_S \frac{f(s)-f(t)}{s-t}ds\pm c\int_S[f(s)-f(t)]
  \delta(s-t)ds=\psi(t)\,\,\,a.e.,
\ee
 because 
 \be\label{e14}  
 \int_S[f(s)-f(t)] \delta(s-t)ds=0\,\,\,a.e.
\ee
 Let us explain \eqref{e14}. 
  The union of the discontinuity points of  a summable function $f$ and the points where the function $f$ is not defined   has measure zero. The 
  $L^1(S)$ norm of a summable function   is not changed if this function is changed on the set of measure zero. This yields the conclusion that  equation \eqref{e14} holds.
  
  Usually the distribution $\delta(s-t)$ is defined as a linear bounded functional on the space of continuous functions with the $\sup$-norm.
  On the other hand, one can define $\delta(s-t)$ as the kernel of the
  identity operator $I$  in the space $L^p(S), \, p\ge 1$. 
    In this case,  the relation $\int_S f(s)\delta(s-t)ds=f(t)$ holds $a.e.$  The expression  $f(t)$ for $f\in L^p(S)$, $p\ge 1$, is not defined at every point $t$, but is  defined almost everywhere. The set of points at which $f$
    is discontinuous, or is not defined (including the points at which $f=\infty$)
    is of measure zero. Outside of this set, that is, a.e., one has
   \be\label{e14a}   I f=\int_{S} \delta(s-t) f(s)ds=f(t) \quad a.e.,
   \ee
 where $I$ is the identity operator with the kernel $\delta(s-t)$.
 Formula \eqref{e14a} gives another way to see  that formula \eqref{e14} is valid. 
 Since in the general theory of distributions the product of two distributions
 is not defined (including the case when one of the distributions is a 
 function,  belonging to $L^p(S)$, $p\ge 1$), formula  \eqref{e14a}
 is useful for $f\in L^p(S)$, $p\ge 1$.
 
  Finally, recall  that for a continuous function $f$ on $S$ 
 equation \eqref{e14} is obvious. If $f\in L^1(S)$ then  the Luzin's theorem,
 see, for example, \cite{S}, p.157,  says that for any
 $f\in L^1(S)$, where $S$ is bounded, and for an arbitrary small $\epsilon>0$, there exists a closed subset $Q$ of $S$,
  $meas Q>meas S-\epsilon$, on which $f$ is continuous. Since $\epsilon>0$
  is arbitrary small, equation \eqref{e14} holds.
 
  Lemma 1 is proved. \hfill$\Box$

 {\bf Remark 4.} Let $S=S_0$ be the unit circle on the complex plane and
 \be\label{e15}  
  B_0f:=\frac 1 {i\pi}\int_{S_0}\frac{f(s)}{s-t}ds.
\ee
For an arbitrary $f\in L^1(S_0)$ the $B_0f$ may be more singular than $f$.

{\bf Example 1.} Let 
\be\label{e16}  
  f:=\sum_{n=1}^\infty \frac{e^{in\phi}-e^{-in\phi}}{n}=2i\sum_{n=1}^\infty \frac{\sin n\phi}{n}=2i \frac{\pi-\phi}2,\quad 0<\phi<2\pi,
\ee
see \cite{GR}, p.52.
 The function
\eqref{e16} is a bounded function on the interval $0\le \phi \le 2\pi$.
 
 One can check,  that
  \be\label{e17} 
  B_0f:=\sum_{n=1}^\infty \frac{B_0e^{in\phi}-B_0e^{-in\phi}}{n}=
  \sum_{n=1}^\infty \frac{e^{in\phi}+e^{-in\phi}}{n}=2\sum_{n=1}^\infty \frac{\cos n\phi}{n}=-\frac 1 2 \ln(2-2\cos \phi),\quad 0<\phi<2\pi.
\ee
see \cite{GR}, p.52. 
The validity of the change of the order of application of $B_0$ and the summation with respect to $n$
follows from the convergence of the series $\sum_{n=1}^\infty \frac{\cos n\phi}{n}$.

 We have used the formulas,  see \cite{MP}, p. 47,
\be\label{e18} 
  B _0e^{in\phi}=e^{in\phi},\quad n>0;   \quad B_0e^{-in\phi}=-e^{-in\phi},\quad n>0.
\ee
Thus,  $B_0f$ has a singularity at $\phi=0$ on the interval 
$0\le \phi \le 2\pi$, while $f$ is a bounded function. 

{\bf Remark 5.} One can check that if $f\in L^1(\R)$ and $Sf:=\frac 1 {2\pi i}\int_{\R}\frac{f(s) ds}{s-t}$, then for some $f$ the $Sf$ is more singular than $L^1(\R)$.
Indeed, consider $Sf$ as a convolution in the sense of distributions. 

Take the Fourier transform of $Sf$ and use the formula $F (f* \frac1s)=Ff( i\pi {sign} \xi)$, where $Ff=\int_{\R}e^{i\xi s}f(s)ds$, $*$
denotes the convolution: $f*h:= \int_{\R}f(s)h(t-s)ds$, $Ff*h=F(f) F(h)$. One has $F\frac 1 s=i\pi{ sign }\xi$,  see \cite{GS}. Since the set $Ff$ is the set
of continuous functions of $\xi$ when $f$ runs through $L^1(\R)$, the set
of functions $\{Ff \times  {sign}\xi\}$ is {\em not} the set of continuous functions of $\xi$. Therefore, $Sf$ does not belong to $L^1(\R)$ for some $f\in L^1(\R)$.

  \section{Conclusion}\label{S:2}
  Let $D$ be a connected bounded domain in $\R^2$, $S$ be its  boundary, which is closed, connected and smooth.  Let $\Phi=\Phi(z),\,\,z=x+iy,$ be an analytic function in $D$. 

 A necessary and sufficient condition for a function $f=f(s)\in L^1(S)$  to be a boundary value of an analytic in $D$ function $\Phi(z)$ is  equation \eqref{e11}.

A necessary and sufficient condition for a function $f=f(s)\in L^1(S)$  to be a boundary value of an analytic in $D'$ function $\Phi(z)$, $\Phi(\infty)=0$, is  equation \eqref{e11a}.


\end{document}